\documentclass{article}

\usepackage{arxiv}

\usepackage[T1]{fontenc}    
\usepackage{doi}

\usepackage{natbib}

\usepackage{graphicx,amsmath}
\usepackage{tikz}
\usepackage{booktabs}
\usepackage{colortbl}
\usepackage{xcolor}
 
\newcommand*{\MinNumber}{0}%
\newcommand*{\MaxNumber}{30}%

\newcommand{\ApplyGradient}[1]{%
        \pgfmathsetmacro{\PercentColor}{100.0*(#1-\MinNumber)/(\MaxNumber-\MinNumber)}
        \edef\temp{\noexpand\cellcolor{gray!\PercentColor!white}}\temp#1
       }

\newcommand{\ApplyGradientMaxOutlier}[1]{%
        \setlength{\fboxsep}{0pt}
        \pgfmathsetmacro{\PercentColor}{100.0*(\MaxNumber-\MinNumber)/(\MaxNumber-\MinNumber)}
        \edef\temp{\noexpand\cellcolor{gray!\PercentColor!white}}\temp#1
}

\newcolumntype{L}[1]{>{\raggedright\let\newline\\\arraybackslash\hspace{0pt}}m{#1}}
\newcolumntype{C}[1]{>{\centering\let\newline\\\arraybackslash\hspace{0pt}}m{#1}}
\newcolumntype{R}[1]{>{\raggedleft\let\newline\\\arraybackslash\hspace{0pt}}m{#1}}

\title{Order acceptance and scheduling in capacitated job shops}

\newif\ifuniqueAffiliation

\ifuniqueAffiliation 
\author{ \href{https://orcid.org/0000-0000-0000-0000}{\includegraphics[scale=0.06]{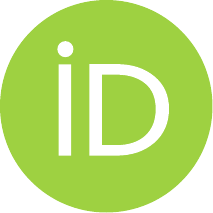}\hspace{1mm}David S.~Hippocampus}\thanks{Use footnote for providing further
		information about author (webpage, alternative
		address)---\emph{not} for acknowledging funding agencies.} \\
	Department of Computer Science\\
	Cranberry-Lemon University\\
	Pittsburgh, PA 15213 \\
	\texttt{hippo@cs.cranberry-lemon.edu} \\
	\And
	\href{https://orcid.org/0000-0000-0000-0000}{\includegraphics[scale=0.06]{orcid.pdf}\hspace{1mm}Elias D.~Striatum} \\
	Department of Electrical Engineering\\
	Mount-Sheikh University\\
	Santa Narimana, Levand \\
	\texttt{stariate@ee.mount-sheikh.edu} \\
}
\else
\usepackage{authblk}

\newbox{\orcid}\sbox{\orcid}{\includegraphics[scale=0.06]{orcid.pdf}} 
\author[1]{%
	\href{https://orcid.org/0000-0002-5869-0425}{\usebox{\orcid}\hspace{1mm}Florian Linß\thanks{\texttt{florian.linss@tu-dresden.de}}}%
}
\author[2]{%
	\href{https://orcid.org/0000-0002-9786-677X}{\usebox{\orcid}\hspace{1mm}Mike Hewitt\thanks{\texttt{mhewitt3@luc.edu}}}%
}
\author[1]{%
	\href{https://orcid.org/0000-0003-0753-0517}{\usebox{\orcid}\hspace{1mm}Janis S. Neufeld\thanks{\texttt{janis\_sebastian.neufeld@tu-dresden.de}}}%
}
\author[1]{%
	\href{https://orcid.org/0000-0003-4711-2184}{\usebox{\orcid}\hspace{1mm}Udo Buscher\thanks{\texttt{udo.buscher@tu-dresden.de}}}%
}
\affil[1]{TU Dresden}
\affil[2]{Loyola University Chicago}
\fi

\renewcommand{\shorttitle}{Order acceptance and scheduling in capacitated job shops}

\hypersetup{
pdftitle={Order acceptance and scheduling in capacitated job shops},
pdfsubject={cs.DM, math.OC},
pdfauthor={Florian Linß, Mike Hewitt, Janis S. Neufeld, Udo Buscher},
pdfkeywords={scheduling; capacity planning; OR in agriculture},
}

\begin{document}
\maketitle

\begin{abstract}
We consider a capacitated job shop problem with order acceptance. This research is motivated by the management of a research and development project pipeline for a company in the agricultural industry whose success depends on regularly releasing new and innovative products. The setting requires the consideration of multiple problem characteristics not commonly considered in scheduling research. Each job has a given release and due date and requires the execution of an individual sequence of operations on different machines (job shop). There is a set of machines of fixed capacity, each of which can process multiple operations simultaneously. Given that typically only a small percentage of jobs yield a commercially viable product, the number of potential jobs to schedule is in the order of several thousands. Due to limited capacity, not all jobs can be started. Instead, the objective is to maximize the throughput. Namely, to start as many jobs as possible. We present a Mixed Integer Programming (MIP) formulation of this problem and study how resource capacity and the option to delay jobs can impact research and development throughput. We show that the MIP formulation can prove optimality even for very large instances with less restrictive capacity constraints, while instances with a tight capacity are more challenging to solve.
\end{abstract}

\keywords{scheduling  \and capacity planning \and OR in agriculture}

\section{Introduction}
Job shop scheduling has attracted significant attention in research and practice \cite{Xiong2022}. However, practical planning problems often necessitate integrating further or new requirements into scheduling models to enable real-world applications. Based on a scheduling task from the agricultural industry, we study a job shop scheduling problem in which each resource (or machine) has a given capacity and can therefore process several jobs simultaneously. The capacity of a machine must not be exceeded at any time. In the context of agricultural production, a resource can be a field, a greenhouse, etc.
Following \cite{Nuijten1996}, we refer to this as capacitated job shop scheduling problem. 

Since the jobs arise from research and development and only a small percentage of jobs results in commercially successful products, the total number of jobs to schedule can reach several thousand. At the same time, the available capacity may not be sufficient to process all jobs. This leads to an order acceptance and scheduling problem \cite{Slotnick2011}. The objective of order acceptance is to maximize throughput, i.e. the number of orders started in the planning period under consideration. 

\cite{Nuijten1996} are the first to introduce a formulation for the capacitated job store problem and solve it as a constraint satisfaction problem. However, there is little other literature on this problem. \cite{Verhoeven1998} propose a similar problem as a resource-constrained scheduling problem (RCSP). In their model, resources can process several operations simultaneously. At the same time, operations might require multiple resources for processing, and partly alternative resources are available. However, the order acceptance decision is not considered. Several other RCSP studies include aspects of the problem at hand. For a recent overview on RCSP see \cite{hartmann2022}.

The order acceptance and scheduling problem in job shop environments is discussed by \cite{Ebben2005}. However, machines are not able to process several jobs at a time. \cite{Wang2021} and \cite{Christ2022} also consider order acceptance in a job shop but on a more aggregated level, i.e., jobs are only assigned to certain planning periods without detailed sequencing.
\cite{Lei2017} propose a neighborhood search for order acceptance in an (uncapacitated) job shop. To the best of our knowledge, this combined order acceptance and capacitated job shop scheduling problem has not been studied in the literature until now and existing planning approaches cannot be applied directly.

In this study, we propose an efficient time-indexed mixed-integer programming model (MIP). In a computational study, we show that the MIP formulation can solve even large instances. Hence, it can already provide decision support in a real-world setting. We also analyze the influence of several factors, such as the capacity of the machines and the tightness of due dates. The gained insights allow a better understanding of the problem.

The remainder is structured as follows: Section \ref{sec:problem} defines the studied problem and presents a mixed-integer programming model. We analyze the MIP by solving several instances of various sizes in a computational study in Section \ref{sec:study}. Finally, we discuss the results and point to directions for future research in Section \ref{sec:conclusion}.

\section{Problem definition and MIP formulation}\label{sec:problem}
The used notation is displayed in Table \ref{tab1}. We consider the scheduling of $n$ jobs, with job $j \in J$. Completing a job requires the completion of a sequence of operations on $m$ machines, with machine $i \in M$. The sequence of machines may differ for each job (i.e., job shop) and is given by $(\sigma_{1j}, \sigma_{2j},...,\sigma_{mj})$. Completing each operation of job $j$ requires the usage of a portion of a machine's $i$ capacity $q_{ij}$ for a fixed processing time $p_{ij}$. This means that different jobs may occupy a machine simultaneously up to a given capacity per machine of $Q_{i}$. Each job $j$ has a release date $r_j$ at which it can be started and a due date $d_j$ by which its last operation must be completed. We note that given the release and due dates and the processing times, we can derive a time window during which each operation of a job can be started. We denote such a time window by $T_{\sigma_{ij}}=[\alpha_{\sigma_{ij}},\beta_{\sigma_{ij}}]$, where
$\alpha_{\sigma_{ij}} = r_j + \sum^{i-1}_{i'=1}p_{\sigma_{i'j},j}$ and $\beta_{\sigma_{ij}} = d_j - \sum^m_{i'=i}p_{\sigma_{i'j},j}$.
With this, there is a fixed planning horizon $H$, resp. $T=\{1,\ldots,H\}$. 

\begin{table}\centering\footnotesize
\caption{Notation}\label{tab1}
\begin{tabular}{ll}
\hline
$n$ & Total number of jobs $j \in J$\\
$m$ & Total number of machines $i \in M$ \\
$r_j$ & Release date of job $j \in J$ \\
$d_j$ & Due date of job $j \in J$ \\
$\sigma_{ij}$ & Machine required for operation $i$ of job $j$ \\
$p_{ij}$ & Processing time of job $j$ on machine $i$\\
$q_{ij}$ & Capacity usage of job $j$ on machine $i$\\
$Q_{i}$ & Available capacity of machine $i$\\
$\alpha_{\sigma_{ij}}$ & Beginning of time window for operation $i$ of job $j$\\
$\beta_{\sigma_{ij}}$ & End of time window for operation $i$ of job $j$\\
$x_{ijt}$ & Binary variable; 1 if job $j$ starts at time $t$ on machine $i$\\
$z_{j}$ & Binary variable; 1 if job $j$ is started / accepted\\
\hline
\end{tabular}
\end{table}

As there may be insufficient capacity to start and complete all jobs, the primary objective is the number of jobs started, which is to be maximized. We refer to this objective as the \textit{throughput} objective.

To define the mathematical model, we let the binary variable $z_{j}, j \in J$ denote whether job $j$ is begun. Let the time-indexed binary variable $x_{ijt}, i \in M, j \in J, t \in T$, denote whether job $j$ is begun on machine $i$ in date $t$. Adapted from \cite{Ku2016}, the linear mixed-integer programming model is defined as follows: 

\begin{equation}
    \text{maximize}  \quad  \sum_{j \in J} z_{j}
\end{equation}
subject to
\begin{equation}
\sum_{t\in T_{\sigma_{ij}}} x_{ijt} = z_{j} \;\;\; \forall j \in J, i=1,\ldots,m
\label{eq:ensure_each_step_start}
\end{equation}
\begin{equation}
\sum_{t\in T_{\sigma_{i-1,j}}} (t+ p_{\sigma_{i-1,j},j})\cdot x_{\sigma_{i-1,j},jt}  \leq \sum_{t'\in T_{\sigma_{ij}}} t'\cdot x_{\sigma_{ij},jt'} \;\;\; \forall j \in J, i=2,\ldots,m
\label{eq:ensure_step_schedule}
\end{equation}
\begin{equation}
\sum_{j \in J}
\sum_{\tau \in \Theta_{ijt}} q_{ij}\cdot x_{ij\tau} \leq Q_{i} \;\;\; \forall i \in M, t\in T, \text{where } \Theta_{ijt} = \{t-p_{ij}+1,\ldots,t\}
\label{eq:ensure_resource_capacity}
\end{equation}
\begin{equation}
z_{j} \in \{0,1\}, \forall j \in J
\label{eq:define_z_vars}
\end{equation}
\begin{equation}
x_{ijt} \in \{0,1\}, \forall j \in J, i\in M, t\in T_{\sigma_{ij}}
\label{eq:define_d_vars}
\end{equation}
Constraints (\ref{eq:ensure_each_step_start}) ensure that each of its operations is started if a job is started. Constraints (\ref{eq:ensure_step_schedule}) ensure an operation of a job is not started until after the preceding operation is completed. Constraints (\ref{eq:ensure_resource_capacity}) ensure the capacity of a resource used by jobs is not exceeded at any date. Constraints (\ref{eq:define_z_vars}) - (\ref{eq:define_d_vars}) define the variables and their domains. 
\\

\section{Computational study}\label{sec:study}
The proposed model was implemented and run on an Intel(R) Xenon(R) CPU E5-4627 with 3.3 GHz clock speed and 768 GB RAM using CPLEX 20.1.0 with max. 4 parallel threads. Runtime was limited to 1,200 seconds. Based on the real-world problem, we generated 1,680 instances with $m=5$ machines and $n \in \{30,50,100,250,500,1000,2000\}$ jobs. All parameters have integer values. The processing times $p_{ij}$ were taken from a uniform distribution within the range $[1,5]$, while the capacity usage $q_{ij}$ is in the range $[20,30]$, and the release dates $r_j$ from $[1,20]$. The due dates are determined by $d_j=r_j+\sum_{\forall i} p_{ij} + w$, where $w$ is a parameter to define the length of the time window, $w\in \{10,20,30\}$. We introduce a capacity factor $f\in \{0.7,0.8,0.9,1.0,1.2,1.5,2.0,5.0\}$ to analyze the impact of the machines' capacities. The machines' capacities are then determined for each instance by multiplying $f$ with the average total required capacity usage for all jobs over the planning horizon. The capacity factor $f$ can be interpreted as the expected acceptance rate if all jobs were equal and evenly distributed over the planning horizon. In addition, we generated instances in which each machine can only process a single job at a time. This results in a classical job shop problem.
\begin{table}\centering\footnotesize
\caption{Average results depending on the number of jobs}\label{tab:Results}
\begin{tabular}{rrrrr}
\toprule
\#Jobs & Acc. Rate [\%] & Gap [\%] & Runtime [s] & \# optimal sol. \\ \midrule
30 & 76.8 & 1.0 & 255.7 & 211 / 240 \\
50 & 80.8 & 5.4 & 709.0 & 109 / 240\\
100 & 84.9 & 7.4 & 749.0 & 91 / 240 \\
250 & 86.5 & 5.0 & 755.9 & 98 / 240 \\
500 & 87.2 & 3.8 & 818.9 & 89 / 240 \\
1,000 & 83.8 & 12.8 & 931.3 & 71 / 240 \\
2,000 & 76.1 & 42.0 & 1,063.9 & 47 / 240 \\ \midrule
Total & 82.3 & 11.0 & 754.8 & 716 / 1,680 \\ \bottomrule
\end{tabular}
\end{table}

Table \ref{tab:Results} shows the average results depending on the number of jobs. For instances with up to 500 jobs, the average MIP gap is less than 7.5\%; for $n=500$ it is even 3.8\%. For $n=$2,000, 47 out of 240 instances could still be solved optimally. This shows that the MIP formulation is able to solve very large instances.

\begin{table}\centering\footnotesize
\caption{Influence of the capacity factor on the average optimality gap depending on the number of jobs and average acceptance rate per capacity factor}\label{tab:CapFactor}
\begin{tabular}{rR{0.1cm}R{1.3cm}R{0.9cm}R{0.9cm}R{0.9cm}R{0.9cm}R{0.9cm}R{1.0cm}R{0.9cm}R{0.9cm}}
\toprule
Cap. $f$ & & Job shop & 0.7 & 0.8 & 0.9 & 1.0 & 1.2 & 1.5 & 2.0 & 5.0 \\ \midrule
30 & & \ApplyGradient{0.2} & \ApplyGradient{0.0} & \ApplyGradient{0.0} & \ApplyGradient{0.0} & \ApplyGradient{2.2} & \ApplyGradient{5.0} & \ApplyGradient{0.8} & \ApplyGradient{0.0} & \ApplyGradient{0.0}
\\50 & & \ApplyGradient{0.0} & \ApplyGradient{6.0} & \ApplyGradient{9.9} & \ApplyGradient{9.5} & \ApplyGradient{8.7} & \ApplyGradient{9.0} & \ApplyGradient{0.1} & \ApplyGradient{0.0} & \ApplyGradient{0.0}
\\100 & & \ApplyGradient{7.0} & \ApplyGradient{14.9} & \ApplyGradient{14.0} & \ApplyGradient{13.7} & \ApplyGradient{12.0} & \ApplyGradient{4.2} & \ApplyGradient{0.0} & \ApplyGradient{0.0} & \ApplyGradient{0.0}
\\\#Jobs\quad250 & & \ApplyGradient{12.9} & \ApplyGradient{11.6} & \ApplyGradient{9.6} & \ApplyGradient{8.8} & \ApplyGradient{7.9} & \ApplyGradient{1.8} & \ApplyGradient{0.0} & \ApplyGradient{0.0} & \ApplyGradient{0.0}
\\500 & & \ApplyGradient{16.9} & \ApplyGradient{9.2} & \ApplyGradient{7.1} & \ApplyGradient{5.4} & \ApplyGradient{4.8} & \ApplyGradient{3.2} & \ApplyGradient{0.4} & \ApplyGradient{0.0} & \ApplyGradient{0.0}
\\1,000 & & \ApplyGradient{16.6} & \ApplyGradient{10.5} & \ApplyGradient{5.4} & \ApplyGradient{4.2} & \ApplyGradient{3.2} & \ApplyGradient{2.8} & \ApplyGradientMaxOutlier{76.3} & \ApplyGradient{0.0} & \ApplyGradient{0.0}
\\2,000 & & \ApplyGradient{29.9} & \ApplyGradient{22.1} & \ApplyGradient{20.9} & \ApplyGradient{17.0} & \ApplyGradient{14.3} & \ApplyGradient{14.9} & \ApplyGradientMaxOutlier{223.1} & \ApplyGradient{23.8} & \ApplyGradient{0.0} \\ \midrule
Total & & 11.9 & 10.6 & 9.6 & 8.4 & 7.6 & 5.8 & 43.0 & 3.4 & 0.0 \\
\midrule \midrule
Acc. Rate & & 19.0 & 58.6 & 66.2 & 73.7 & 80.2 & 91.1 & 90.6 & 97.8 & 100.0  \\ \bottomrule
\end{tabular}
\end{table}

Table \ref{tab:CapFactor} displays the average optimality gap of CPLEX and the acceptance rates dependent on the capacity factor $f$. The job shop case, where no jobs can be processed in parallel by each machine, shows the highest gap values of, on average, 11.9\%. Using different capacity factors $f$ influences the solution process. Especially medium and large problems show decreased gap values for increasing capacity factors. This indicates that the less tight the capacities are, the easier it is to solve the problem. For example, all instances with up to 1,000 jobs are solved optimally for $f=2.0$. For very small instances ($n\leq 100$), another influencing factor on the gap values can be observed: if the number of jobs is low, the acceptance of a single job has a relatively high impact on the objective value (e.g., for $n=30$, one job equals $1/30=3.3\%$ of the throughput). This explains the decreasing gap values between small and medium instances (e.g., $8.7\%$ for $n=50$ and $f=1.2$ and $3.2\%$ for $n=1000$).
In contrast, the absolute gap value as the difference between the best-known solution and the upper bound of the number of accepted jobs is comparable to the larger instances. The high gap values for $f=1.5$ and $n\geq 1,000$ are mainly caused by outliers. For these instances, CPLEX can only find weak feasible solutions, which leads to extreme gap values. 

The acceptance rate in the last row of Table \ref{tab:CapFactor} is lower than the capacity factor $f$. On average, only 19\% of all jobs can be accepted in the job shop case, while a capacity factor of 5.0 is necessary to accept all jobs for all instances. Interestingly, the size of the time window $w$ does not have a relevant impact on the acceptance rate of the problem.

\section{Conclusion and future research}\label{sec:conclusion}
We present a new order acceptance and capacitated job shop scheduling problem motivated by a real-world planning problem in the agricultural industry. Several instances of the proposed MIP formulation with up to 2,000 jobs could be solved optimally. In many practical instances, the proposed formulation can provide decision support. Within the computational study, we could see that especially the relation between the number of jobs and the available capacity (i.e., the capacity factor or the acceptance rate, respectively) strongly influences the complexity of the problem. While instances where (nearly) all jobs can be accepted are relatively easy to solve, it is more difficult to find the optimal solution for problems with low capacity.

Besides throughput, other objectives are relevant in practice and should be addressed in future research. Maximizing the total profit with differing profits per job can reflect the importance of the corresponding project. Another aspect of the application in agricultural business is that waiting times between operations should be minimized. From a methodological point of view, constraint programming has shown good performance for several scheduling problems \cite{Naderi2023} and could be tested for this problem, too. To solve very large instances arising in practice, which can have tens of thousands of jobs, the development of heuristic solution approaches is worthwhile. 

\bibliographystyle{apalike} %

\end{document}